\documentclass[10pt,a4paper]{amsart}
\usepackage{amsfonts,amsmath,amssymb}
\usepackage{hyperref}
\usepackage[all]{xy}
\newtheorem*{theorem*}{Theorem}
\newtheorem{lemma}{Lemma}[section]
\newtheorem{proposition}[lemma]{Proposition}
\newtheorem{remark}[lemma]{Remark}

\newtheorem{theorem}[lemma]{Theorem}
\newtheorem{definition}[lemma]{Definition}
\newtheorem{notation}[lemma]{Notation}

\newtheorem{corollary}[lemma]{Corollary}

\newtheorem*{conjecture*}{Conjecture}
\baselineskip 18pt \textwidth 16cm  \theoremstyle{plain}

\newcommand{\tr}{\operatorname{Tr}}

\newcommand{\swrtz}{\mathcal{S}}

\newcommand{\Mat}{\operatorname{Mat}}

\newcommand{\eps}{\varepsilon}

\newcommand{\id}{\operatorname{Id}}
\newcommand{\re}{\operatorname{Re}}
\renewcommand{\Im}{\operatorname{Im}}

\newcommand{\R}{{\mathbb R}}
\newcommand{\C}{{\mathbb C}}

\newcommand{\Fre}{{Fr\'{e}chet \,}}

\newcommand{\g}{{\mathfrak{g}}}

\newcommand{\GL}{\operatorname{GL}}

\newcommand{\gl}{{\mathfrak{gl}}}

\newcommand{\Sc}{{\mathcal S}}

\begin{document}

\title[Uniqueness of Shalika functionals]{Uniqueness of Shalika functionals (the Archimedean case)}
\date{\today}
\author{Avraham Aizenbud}
\address{Avraham Aizenbud, Faculty of Mathematics
and Computer Science, The Weizmann Institute of Science POB 26,
Rehovot 76100, ISRAEL.} \email{aizenr@yahoo.com}
\author{Dmitry Gourevitch} \address{Dmitry Gourevitch, Faculty of Mathematics
and Computer Science, The Weizmann Institute of Science POB 26,
Rehovot 76100, ISRAEL.} \email{guredim@yahoo.com}
\author{Herv{\'e} Jacquet} \address{Herv{\'e} Jacquet, Mathematics
Department of Columbia University, MC 4406, 2990 Broadway New
York, NY 10027} \email{hj@math.columbia.edu}
%
%
%
%
%
%
%
%
%
%
\maketitle

\begin{abstract}

Let $F$ be either $\R$ or $\C$. Let $(\pi,V)$ be an irreducible
admissible smooth \Fre representation of $GL_{2n}(F)$.  A Shalika
functional $\phi:V \to \C$ is a continuous linear functional such
that for any $g\in GL_n(F), \, A \in \Mat_{n \times n}(F)$ and $v\in
V$ we have $$ \phi\left[\pi\begin{pmatrix}
  g & A \\
  0 & g
\end{pmatrix})v\right] = \exp(2\pi i \re( \tr (g^{-1}A))) \phi(v).$$
In this paper we prove that the space of Shalika functionals on
$V$ is at most one dimensional.

For non-Archimedean $F$ (of characteristic zero) this theorem was
proven in \cite{JR}.
\end{abstract}

\setcounter{tocdepth}{2}
 \tableofcontents

\section{Introduction}
Let $F$ be either $\R$ or $\C$. Let $(\pi,V)$ be an admissible
smooth \Fre representation of $GL_{2n}(F)$. We assume that $V$ is
the canonical completion of an irreducible Harish-Chandra $(\g,K)$-
module in the sense of Casselman-Wallach (see e.g. \cite{Wal2},
chapter 11). A \textbf{ Shalika functional} $\phi:V \to \C$ is a
continuous linear functional such that for any $g\in GL_n(F), \, A
\in Mat_{n \times n}(F)$ and $v\in V$ we have $$
\phi\left[\pi\begin{pmatrix}
  g & A \\
  0 & g
\end{pmatrix}v\right] = \exp\left(2\pi i \re( \tr (g^{-1}A))\right) \phi(v).$$
In this paper we prove the following theorem.

\begin{theorem}\label{UShal}
 Let $(\pi,V)$ be an irreducible admissible
smooth \Fre representation of $GL_{2n}(F)$. Then the space of
Shalika functionals on $V$ is at most one dimensional.
\end{theorem}

For non-Archimedean $F$ (of characteristic zero) this theorem was
proven in \cite{JR}. The proof in \cite{JR} is based on the fact
that $(\GL_{2n}(F),\GL_{n}(F) \times \GL_{n}(F))$ is a Gelfand
pair, which was also proven in \cite{JR}, and the method of
\cite[Section 3]{FJ} of integration of Shalika functionals.

In the Archimedean case those two ingredients also exist. Namely,
\cite[Section 3]{FJ} is valid also in the Archimedean case, and
the fact that $(\GL_{2n}(F),\GL_{n}(F) \times \GL_{n}(F))$ is a
Gelfand pair is proven in \cite{AG_RJR}.

The proof that we present here is similar to the proof in \cite{JR}.
The main difierence is that we have to prove the continuity of a
certain linear form.

\subsection{Structure of the proof}$ $\\
We construct a linear map from the space of Shalika functionals to
the space of linear periods (linear functionals on $V$ that are
invariant by $\GL_{n}(F) \times \GL_{n}(F)$)  and prove that the map
is injective. Hence the uniqueness of the linear periods implies
uniqueness of the Shalika functionals. The uniqueness of linear
periods, i.e. the fact that $(\GL_{2n}(F),\GL_{n}(F) \times
\GL_{n}(F))$ is a Gelfand pair, is proven in \cite{AG_RJR}.

\subsection{Structure of the paper} $ $\\
In \S \ref{Prel} we fix notation and terminology. In \S
\ref{SecIntShal} we describe a way of obtaining a linear period from
a Shalika functional by integration, as in \cite[Section 3]{FJ}. In
\S \ref{SecProp} we investigate the properties of the obtained
period. In \S \ref{SecUShal} we explain how this implies the
uniqueness of Shalika functionals.

\subsection{Acknowledgements} $ $\\
Aizenbud and Gourevitch thank {\bf Josef Bernstein}, {\bf Wee Teck
Gan} and {\bf Binyong Sun} for useful remarks.

Aizenbud and Gourevitch were partially supported by a BSF grant, a
GIF grant, and an ISF Center of excellency grant.
Aizenbud was also partially supported by ISF grant No. 583/09.

\section{Preliminaries and notation} \label{Prel}
\subsection{Notation} \label{Conv}
\begin{itemize}
\item Henceforth we fix an Archimedean field $F$ (i.e. $F$ is $\R$ or $\C$).
\item  For a group $G$ acting on a vector space $V$ we denote by $V^{G}$ the space of $G$-invariant
vectors in $V$. For a character $\chi$ of $G$ we denote by
$V^{G,\chi}$ the space of $(G,\chi)$-equivariant vectors in $V$.
\item For a smooth real algebraic variety $M$ we denote by
$\Sc(M)$ the space of Schwartz functions on $M$, i.e. the space of
smooth functions that are rapidly decreasing as well as all their
derivatives. For precise definition see e.g. \cite{AG_Sc}.

\item We fix a natural number $n$ and denote $G:=GL_{2n}(F)$.
\item We fix a norm on $G$ by $$||g||:= \sum _{1\leq i,j \leq 2n} |g_{ij}|^2+ \sum _{1\leq i,j \leq
2n}|(g^{-1})_{ij}|^2.$$
\item We denote
$$G_1:= \left\{\begin{pmatrix}
  g & 0 \\
  0 & Id
\end{pmatrix} \,|\, g \in GL_n(F)\right \} \subset G $$

\item We denote by $\nu:GL_n(F) \to G_1$ the isomorphism
defined by
$$\nu(g):= \begin{pmatrix}
  g & 0 \\
  0 & Id
\end{pmatrix}.$$
Note that for any $X \in \Mat(n \times n,F)$, $d\nu(X)=
\begin{pmatrix}
  X & 0 \\
  0 & 0
\end{pmatrix}.$

\item We denote $$H:= \left\{\begin{pmatrix}
  g & 0 \\
  0 & g
\end{pmatrix} \,|\, g \in GL_n(F)\right \} \subset G$$

\item We denote $$U:= \left\{\begin{pmatrix}
  Id & A \\
  0 & Id
\end{pmatrix} \,|\, A \in Mat_{n\times n}(F)\right \}\subset G $$
\item We denote by $\mu:\Mat(n \times n,F) \to U$ the isomorphism
defined by
$$\mu(A):= \begin{pmatrix}
  Id & A \\
  0 & Id
\end{pmatrix}.$$
Note that for any $X \in \Mat(n \times n,F)$, $d\mu(X)=
\begin{pmatrix}
  0 & X \\
  0 & 0
\end{pmatrix}.$

\item We denote by $\tau : U \to F$ the homomorphism given  by
$$\tau (\mu(A)) := \tr(A).$$

\item We let $\psi$ be the additive character of $F$ defined by $\psi(x):=e^{2\pi \mathrm{i}\re x}$. We
define an homomorphism $\Psi: U\rightarrow F^\times$ by
$$ \Psi:= \psi \circ \tau.$$
We extend $\Psi$ to an homomorphism $\Psi: HU\rightarrow F^\times$
trivial on $H$.

\item We denote by $K$ the standard maximal compact subgroup of $G$.
Thus $K=O(2n)$ if $F=\R$ and $K= U(2n)$ if $F=\C$.
\end{itemize}

\subsection{Admissible representations}$ $\\
In this paper we consider admissible smooth \Fre representations
 of $G$, i.e. smooth admissible representations $(\pi,V)$ of $G$ such that $V$ is a
 \Fre space and, for any continuous semi-norm $\alpha$ on $V$, there exist another
continuous semi-norm $\beta$ on $V$ and a natural number $M$ such
that for any $g\in G$,
$$\alpha(\pi(g)v) \leq \beta(v)||g||^M.$$

By Casselman - Wallach theorem (see e.g. \cite{Wal2}, chapter 11),
$V$ may be regarded as the canonical model of an irreducible
Harish-Chandra $(\mathfrak{g},K)-$module. By Casselman embedding
theorem (\cite{Cas}), $V$ can be realized as a closed subspace of
a principal series representation. We denote by $\widetilde{V}$
the canonical model of the contragredient Harish-Chandra
$(\mathfrak{g},K)-$module. It is a subspace of the topological
dual $V^*$ of $V$.

\section{Integration of Shalika functionals} \label{SecIntShal}
In this section we fix:
\begin{itemize}
\item an irreducible admissible smooth
\Fre representation $(\pi,V)$ of $G$
\item a Shalika functional $\lambda$ on $V$, i.e. $\lambda \in
(V^*)^{HU,\Psi}$. \end{itemize}

\begin{theorem} \label{IntShal}
There exists $M \in \R$ such that for any $v\in v$ and over the
region of $s\in \C$ with $\re(s)
>M$, the integral
$$L_{\lambda, v}(s):= \int_{g \in G_1}\lambda(\pi(g)v)|\det(g)|^{s-\frac{1}{2}}dg$$
converges absolutely and is a holomorphic function of $s$.

Moreover, $L_{\lambda, v}(s)$ has meromorphic continuation to the
complex plane and is a holomorphic multiple of the $L$-function
$L_{\pi}$ of the representation $\pi$. Finally, for any $\lambda
\neq 0$ there exists $v \in V$ such that $L_{\lambda, v} =
L_{\pi}$.
\end{theorem}

In \cite[Proposition 3.1]{FJ} this theorem is proven under the
following assumption:\\
(*) There exists a continuous semi-norm $\beta$ on $V$ such that
$|\lambda(\pi(g)v)| \leq \beta(v)$ for any $g \in G$.

This may not be true in general. However, we have the following
result.
\begin{lemma} \label{betaG1}
There exist $M >0$ and a continuous semi-norm $\beta$ on $V$ such
that $|\lambda(\pi(g)v)| \leq |\det g|^{-M}\beta(v)$ for any $g \in
G_1$.
\end{lemma}
Before proving the Lemma, we check that, with the help of this
Lemma, the proof of Theorem \ref{IntShal} is still valid. Indeed,
the functions $g\mapsto \lambda(\pi(g)v)$ are bounded in \cite{FJ}
and satisfy a sharper estimate (\cite[Lemma 3.1]{FJ}). Here they
satisfy the following estimates.
\begin{lemma} \label{betaG11}
There is a continuous semi-norm $\gamma$ on  $V$ such that, for any
$v\in V$,
\[ |\lambda (\pi(g)v)| \leq |\det b^{-1}a|^{-M} \gamma(v) \]
for
\[ g = u\left(\begin{array}{cc} a & 0 \\ 0 & b\end{array}\right)k
\]
with $a,b\in GL(n,F)$, $u\in U$, $k\in K$. Furthermore, for any
$v\in V$, there is $\Phi_v\in \swrtz(\Mat(n\times n,F))$ such that
\[ |\lambda (\pi(g)v)| \leq \Phi(b^{-1}a) |\det b^{-1}a|^{-M} ,\]
for $g$ of the above form.
\end{lemma}
\begin{proof}
For the first assertion, we have
\[ \lambda(\pi(g)v) = \Psi(u) \lambda ( \pi(\nu(b^{-1}a))\pi(k)v)
.\] Hence
\[ |\lambda (\pi(g)v)|\leq |\det b^{-1} a| ^{-M} \beta(\pi(k)v) .\]
There is another continuous semi-norm $\gamma$ such that, for all
$k\in K$,
\[ \beta(\pi(k)v)\leq \gamma(v) .\]
The first assertion follows.

For the second assertion, we go through the proof of \cite[Lemma
3.1]{FJ} (which is the above estimate with $M=0$) and arrive at once
at the present estimate.
\end{proof}

The proof of Theorem \ref{IntShal} is still valid. The only
modification is that we need to check that, under our weaker
assumption, two integrals in \cite{FJ} which depend on $s\in \C$,
are still absolutely convergent for $\re s>>0$.

The first integral is integral \cite[45]{FJ}:
\[ \int \lambda(\pi(g)v) \Phi(g) |\det g|^{s+ n -\frac{1}{2}}
d^\times g \] where $\Phi\in\swrtz(\Mat(2n\times 2n,F))$. We write
\[ g = \left(\begin{array}{cc} a & x \\ 0 & b\end{array}\right)k.\]
Then
\[ d^\times g = |\det a|^{-n} d^\times ad^\times b dxdk .\]
By Lemma \ref{betaG11}, the integral of the absolute value is
bounded by
\[ \int |\det a|^{\re s-M-\frac{1}{2}} |\det b| ^{\re s+M + n-\frac{1}{2}}
\left| \Phi\left[\left(\begin{array}{cc} a & x \\ 0 &
b\end{array}\right)k\right]\right|  d^\times a d^\times b dxdk.\]
This does converge absolutely for $\re s>>0$.

The second integral is integral \cite[48]{FJ}. It has the form
\[ \int \lambda\left[ \pi\left(\begin{array}{cc} a & 0\\ 0 & \id
\end{array}\right) \pi(x)v\right] |\det a|^{s-\frac{1}{2}} d^\times a d\mu(x)
\]
where $\mu$ is the measure on $SL(2n,F)$ defined by
\[ \int f(x) d\mu(x) =\]
\[ \int f\left[ \left(\begin{array}{cc} b^{-1} & 0 \\ 0 & \id\end{array}\right)
\left(\begin{array}{cc} \id & u \\ 0 &
\id\end{array}\right)\left(\begin{array}{cc} \id & 0 \\ 0 &
b\end{array}\right)k\right]\Upsilon (u,b^{-1},b;k) \ |\det b|^n
d^\times b du dk \] In this formula $k$ is integrated over
$K'=K\cap SL(2n,F)$ and the function $\Upsilon$ is in $ \swrtz
(\Mat( n\times n,F)^3\times K')$. The integral of the absolute
value of the integrand is bounded by
\[ \int \Phi_v(ab^{-2}) |\det a b^{-2}| ^{-M} \ |\Upsilon|
(u,b^{-1},b;k) |\det b|^n d^{\times} b dudk \ |\det a|^{\re s-
\frac{1}{2}} d^\times a .\] After changing $a$ to $a b^2$, the
integral decomposes into a product:
\[\int \Phi_v(a) |\det a|
^{\re s-M-\frac{1}{2}} d^\times a \times \int |\Upsilon|
(u,b^{-1},b;k) |\det b| ^{n+ 2\re s -1} d^\times bdu dk .\] The
first integral converges for $\re s>>0$. The second integral
converges for all $s$.

It remains to prove Lemma \ref{betaG1}. We will prove the following
more general lemma.

\begin{lemma} \label{ForPol} There exists $M_0>0 $ such that,
for any polynomial $P$ on the real vector space $\Mat(n\times n,F)$,
there exists a continuous semi-norm $\beta_P$ on $V$ such that for
any $g \in \GL_n(F)$ we have
$$| \lambda \left  ( \pi (\nu(g)) v \right )| \leq \beta_P(v)
\frac{1}{|P(g)|} |\det g|^{-M_0}.$$
\end{lemma}
\begin{proof}
We have
$$ \lambda( \pi(\mu(X)) v ) = \psi(\tr X) \lambda(v) \quad \forall X \in
\Mat(n\times n,F).$$

We have then
$$ \lambda (d\pi(d\mu(X))v) = 2 \pi \mathrm{i} \re \tr(X)\lambda(v)
\quad \forall X \in \Mat(n\times n,F)$$ and hence
$$ \lambda    ( \pi (\nu(g)) d\pi(d\mu(X))v   ) = 2 \pi \mathrm{i} \re \tr(gX)\lambda(\pi(\nu(g))v)
\quad \forall X \in \Mat(n\times n,F) \text{ and }g \in \GL_n(F).$$
Similarly, if $Q$ is a polynomial on the real vector space
$\Mat(n\times n,F)$, there is an element $X_Q$ of the enveloping
algebra of $ \gl_{2n}(F)$ such that
$$ \lambda    ( \pi (\nu(g)) d\pi(X_Q)v   ) = Q(g) \lambda( \pi (\nu(g))v)
\quad \forall g \in \GL_n(F).$$
We know that there exist a continuous semi-norm $\beta$ on $V$ and a
natural number $M$ such that $|\lambda(\pi(g)v)| \leq \beta(v)
||g||^M$ for any $g\in G$. Therefore for any $g \in \GL_n(F)$ we
have
$$  |  Q(g) \lambda    ( \pi (\nu(g)) v   )   |  =   | \lambda    ( \pi (\nu(g)) d\pi(\mu(X_Q))v   )   |  \leq
\beta(d\pi(X_Q)v)||\nu(g)||^M.$$
Note that $||\nu(g)||^M=P_0(g)|\det g|^{-2M}$ for a suitable
polynomial $P_0$ on the real vector space $\Mat(n\times n,F)$.
Therefore, we have, with $M_0=2M$,
$$    | \lambda    ( \pi (\nu(g)) v   )  | \leq \beta(d\pi(X_Q)v)
\frac{P_0(g)}{|Q(g)|} |\det g|^{-M_0}.$$
We may take $Q$ of the form $Q=P_0P$ where $P$ is another
polynomial. Since $v \mapsto \beta(d\pi(X_Q)v)$ is a continuous
semi-norm the Lemma follows.
\end{proof}



\section{Properties of $L_{\lambda,v}$} \label{SecProp}
\begin{theorem} \label{Bound}
Let $(\pi,V)$ be an irreducible admissible smooth \Fre
representation of $G$. Fix a Shalika functional $\lambda \in
(V^*)^{HU,\Psi}$ and a vector $v \in V$. Then, for any polynomial
$p$, the product $p(s)L_{\lambda,v}(s)$ is bounded at infinity on
every vertical strip of finite width.
\end{theorem}

In \cite[\S\S 3.3]{FJ} the following statement is proven.

\begin{lemma} \label{FJ3_3}
For $\re(s)$ large enough, $L_{\lambda,v}(s)$ is a finite sum of
functions of the type
$$\mathcal{L}_{u,\xi,\Phi}(s):= \int_{g\in G}
\Phi(g)\xi(\pi(g)u) |\det g|^{s+n-\frac{1}{2}} dg,$$ where $\Phi \in
\Sc(\Mat(2n\times 2n ,F))$, $u \in V$, $\xi \in V^*$.
\end{lemma}
Now Theorem \ref{Bound} follows from the following one.

\begin{theorem} \label{BoundGJ}
Let $(\pi,V)$ be an irreducible admissible smooth \Fre
representation of $G$. Let $\Phi \in \Sc(\Mat(2n\times 2n ,F))$, $u
\in V$ and $\xi \in \widetilde{V}$. Then
$\mathcal{L}_{u,\xi,\Phi}(s)$ has a meromorphic continuation to $\C$
whose product by any polynomial is bounded at infinity on any
vertical strip. The continuation is a holomorphic multiple of
$L_{\pi}(s)=L(s,\pi)$. It satisfies the functional equation
$$
\int \widehat{\Phi}(g)\xi(\pi( {}^tg ^{-1})u)|\det g| ^{1-s+
n-\frac{1}{2}} dg = \gamma(s,\pi,\psi) \mathcal{L}_{u,\xi,\Phi}(s)$$
where
$$ \gamma(s,\pi,\psi) : =
\eps(s,\pi,\psi)\frac{L(1-s,\widetilde{\pi})}{L(s,\pi)}.$$ and
$$ \widehat{\Phi}(X):= \int _{\Mat(2n \times 2n, F)}\Phi(Y) \psi(\mathrm{tr}(
XY^t))
dY.$$ Finally, these assertions remain true if $\xi$ is in $V^{*}$
(topological dual of $V$).
\end{theorem}

This theorem is proven in \cite{GJ} in slightly narrower generality:
the vectors $u$ and $\xi$ are $K-$finite and the function $\Phi$ is
the product of a Gaussian function and a polynomial.  For the
convenience of the reader we indicate how to extend the results of
\cite{GJ}.

We will need the following lemma.

\begin{lemma} \label{GtoT}
Let $T \subset G$ be the torus of diagonal matrices. We will also
regard $T$ as the subset $(F^{\times})^{2n}$ of $F^{2n}$. Let
$\chi:T \to \C^{\times}$ be a multiplicative character.  Let
$(\pi,V)$ be the corresponding representation of principal series
of $G$. Let $v \in V$ and $\xi \in \widetilde{V}$. Let $\Phi$ be a
Schwartz function on $\Mat(2n \times 2n, F)$.

Then there exists a Schwartz function $\phi \in \Sc(F^{2n})$ such
that
$$ \int_{g \in G} \Phi(g)\xi(\pi(g)v) | \det
g|^{s+n-\frac{1}{2}} dg = \int_{t \in T} \phi(t)\chi(t)|\det t|^s
dt$$
for any $s \in \C$ such that the integral on the right converges
absolutely.
\end{lemma}
\begin{proof}
Let $N$  denote the group of upper triangular matrices with unit
diagonal. Let $B=TN$ and $\delta_B$ be the module of the group $B$.
Realize $V$ as the space of smooth functions on $G$ that satisfy
$$f(tg) = \chi(t)f(g)\delta_{B}^{1/2}(t) \text{ and } f(ug)=f(g) \text{ for any }t\in T \text{ and }u\in N.$$
Realize also $\widetilde{V}$ in the corresponding way. Then
 $$\xi(\pi(g)v) = \int _{k\in K}v(kg)\xi(k)dk\,, $$ where $K$ is the standard maximal compact subgroup.  Now
\begin{multline*}
\int_G \Phi(g)\xi(\pi(g)v) | \det g|^{s+\frac{n-1}{2}} dg = \int_G
\int_K \Phi(g) v(kg)\xi(k)| \det g|^{s\frac{n-1}{2}} dgdk = \int_G
\int_K \Phi(k^{-1}g) v(g)\xi(k)| \det g|^{s+\frac{n-1}{2}} dgdk
\end{multline*}
To compute this integral we set
\[g = \left( \begin{array}{cccc} a_1 & u_{1,2} &\cdots & u_{1,2n} \\
                                  0&a_2 & \cdots & u_{2,2n}       \\
                                  \cdots & \cdots & \cdots  \\
                                  0&0&\cdots &a_{2n}
\end{array}\right)k'
\,.\] Then
\[ dg = |a_1|^{1-2n} |a_2|^{2-2n} \cdots |a_{2n-1}|^{-1}\otimes d^\times a_i \otimes du_{i,j}
dk'\] We set
$$\phi(a_1,...,a_{2n}) := \int \Phi\left[ k^{-1}\left( \begin{array}{cccc} a_1 & u_{1,2} &\cdots & u_{1,2n} \\
                                  0&a_2 & \cdots & u_{2,2n}       \\
                                  \cdots & \cdots & \cdots  \\
                                  0&0&\cdots &a_{2n}
\end{array}\right)k'\right]v(k')\xi(k)dkdk' \otimes du_{i,j} \,.$$
Clearly $\phi$ is a Schwartz function on $F^{2n}$ and
$$ \int_{g \in G} \Phi(g)\xi(\pi(g)v) | \det
g|^{s+\frac{n-1}{2}} dg = \int_{t \in T} \phi(t)\chi(t)|\det t|^s
dt$$
for any $s \in \C$ such that the integral on the right converges.
\end{proof}

Now we can prove Theorem \ref{BoundGJ}.

\begin{proof}[Proof of Theorem \ref{BoundGJ}]
The representation $(\pi,V)$ is a sub-representation of a principal
series representation determined by a character $\chi$ of $T$ and
the representation $(\widetilde{\pi},\widetilde{V})$ is then a
quotient of the representation determined by $\chi^{-1}$. For $u\in
V$ and $\xi\in \widetilde{V}$ (or $\xi$ in the principal series
determined by $\chi^{-1}$) we have
$$ \mathcal{L}_{u,\xi,\Phi}(s)= \int_{(a_1,..,a_{2n})
 \in {F^{\times}}^{2n}} \phi(a_1, a_2, . . . , a_{2n})\chi_1(a_1)|a_1|^s .
. . \chi_{2n}(a_{2n})|a_{2n}|^s d^{\times} a_1 . . . d^{\times}
a_{2n}\,.$$ The right hand side extends to a meromorphic function
of $s$ and the product of this function by any polynomial is
bounded at infinity in any vertical strip. Moreover, the function
$\phi$ depends continuously on $\Phi, u\in V, \xi\in
\widetilde{V}$. Therefore the analytic continuation depends
continuously on $\Phi,v\in V, \xi\in \widetilde{V}$. By
continuity, it has the properties stated in the Theorem. To extend
further to the case where $\xi$ is in the topological dual $V^*$
we appeal to the Dixmier Malliavin Lemma (\cite{DM}) applied to
the representation of $SL_{2n}(F)$ on $\Sc (\Mat _{2n\times
2n}(F))$ defined by
\[ g_1 \Phi(X):= \Phi(g_{1}^{-1}X) \]
Thus we may assume $\Phi$ is of the form
\[ \Phi(X) = \int _{SL_{2n}(F)}\Phi_1(g_{1}^{-1} X) f(g_1)dg_1\]
where $f_1$ is a $C^\infty$ function of compact support on
$SL_{2n}(F)$. Then
\[ \mathcal{L}_{u,\xi,\Phi} = \mathcal{L}_{u,\xi_1,\Phi_1} \]
where
\[ \xi_1 (v):= \xi (\pi(f_1)v) \,.\]
Now $\xi_1$ is in $\widetilde{V}$ and our assertion follows.

%
\end{proof}
\begin{remark} The previous result with $\xi \in V^*$ is used
without comment in \cite{FJ}, formula (57). This is why we included
a sketch of the proof.
\end{remark}

%

\begin{theorem} \label{Cont}
There exists $M>0$ such that for any even integer $M'\geq 2$ and
any polynomial $p$ on $\C$, there exists a semi-norm $\beta$ on
$V$ such that $ \left |{{pL_{\lambda,v}}|_{M'+M+\mathrm{i}\R}}^{}
\right|\leq \beta(v)$.
\end{theorem}

First, we will prove the following lemma.

\begin{lemma} \label{LamBet}
There exists $M>0$ such that, for any even integer $M'\geq 2$, there
exists a continuous semi-norm $\beta$ on $V$ such that
$|L_{\lambda,v}(s)|\leq\beta(v)$ for $\re s =M'+M$.
\end{lemma}

\begin{proof}
By Lemma \ref{ForPol} there exists $M_0>0$, and for any polynomial
$P$ on $\Mat(n\times n,F)$, a continuous semi-norm $\beta_P$ such
that
$$ \lambda \left  ( \pi (\nu(g)) v \right ) \leq \beta_P(v)
\frac{1}{|P(g)|} |\det g|^{-M_0}.$$
Let $M:=1/2+n^2+M_0$. Let $M'\geq 2$ be an even integer and let $s
\in \C$ with $\re(s)=M+M'$. Let $$P(X)= |\det X|^{M'}
\prod_{i,j}(1+X_{ij}\overline{X}_{ij}).$$ Let
$$ \beta(v):= \beta_P(v) \int_{X \in\Mat(n \times
n,F)}\frac{dX}{\prod_{i,j}(1+X_{ij}\overline{X}_{ij})}.$$ Now
\begin{multline*}
|L_{\lambda,v}|= \left | \int_{GL_n(F)} \lambda (\pi(\nu(g))v) |\det
g|^{s-1/2}dg \right | \leq \int_{GL_n(F)} \beta_P(v)
\frac{1}{|P(g)|} |\det g|^{-M_0} |\det g|^{n^2 + M'+M_0}dg =
\\  \int_{GL_n(F)} \beta_P(v) \frac{1}{|P(g)|} |\det g|^{n^2 +
M'}dg = \int_{X \in \Mat(n\times n,F)} \beta_P(v) \frac{1}{|P(X)|}
|\det X|^{M'}dX=\beta(v).
\end{multline*}
\end{proof}

\begin{proof}[Proof of Theorem \ref{Cont}]
For any $g \in GL_n(F)$ we have
$$L_{\lambda,\pi(\nu(g))v}(s) = |\det (g)|^{1/2-s}L_{\lambda,v}(s).$$

We can apply this to $g=\exp(tX)$, with $t\in \R$ and $X\in
\Mat(n\times n,F)$. We get
$$L_{\lambda,\pi(\nu(g))v}(s)= |\det( \exp(tX))|^{\frac{1}{2}-s}L_{\lambda,v}(s).$$
Differentiating this identity with respect to $t$ at $t=0$, we get
$$L_{\lambda,d\pi(d\nu(X))v}(s) =(\frac{1}{2}-s)c( X) L_{\lambda,v}(s),$$
where $c(X)=\tr X$ if $F=\R$ and $c(X)=2 \re \tr X$ if $F=\C$.
Similarly, for any polynomial $p$ on $\C$ there exists $X_p$ in the
universal enveloping algebra of $\gl_{2n}(F)$ such that
$$L_{\lambda,d\pi(X_p)v}(s) =p(s) L_{\lambda,v}(s).$$
The theorem follows now from Lemma \ref{LamBet}.
\end{proof}

\begin{notation}
Define another representation $\pi^{\theta}$ on the same space $V$
by $\pi^{\theta}(g):=\pi((g^t)^{-1})$. Recall that $\pi^{\theta}
\cong \widetilde{\pi}$.

For any Shalika functional $\lambda: \pi \to \C$ we define
$\lambda^{\theta}:\pi^{\theta} \to \C$ by
$$\lambda^{\theta}(v):=\lambda \left( \pi \begin{pmatrix}
  0_{nn} & Id_{nn} \\
  -Id_{nn} & 0_{nn}
\end{pmatrix} v\right).$$
It is easy to see that $\lambda^{\theta}$ is a Shalika functional
for the rerpresentation $\pi^\theta$.
\end{notation}

\begin{theorem}[\cite{FJ}, Proposition 3.3] \label{FuncEq}
$$ \gamma(s,\pi,\psi) L_{\lambda,v}^{\pi}(s) = L_{\lambda^{\theta},v}^{\pi^{\theta}}(1-s).$$
\end{theorem}

Using Theorem \ref{Cont} we obtain the following corollary.
\begin{corollary} \label{LineBound}
There exists $N<0$ such that for any odd integer $N'\leq -1$ and
any polynomial $p$ on $\C$ there exists a semi-norm $\beta$ on $V$
such that $|pL_{\lambda,v}|_{N'+N+\mathrm{i}\R}|\leq \beta(v)$.
\end{corollary}

\section{Uniqueness of Shalika functionals} \label{SecUShal}

\begin{theorem} \label{FuncCont}
Let $(\pi,V)$ an irreducible admissible representation of $G$. Let
$\lambda$ be a Shalika functional. Then the functional $L(\lambda):V
\to \C$ defined by
$$L(\lambda)(v):= \frac{L_{\lambda,v}}{ L_\pi}(\frac{1}{2})$$
is continuous.
\end{theorem}
\begin{proof}
By Theorem \ref{Cont} we choose $M>1$ such that for any polynomial
$p$ there exists a semi-norm $\beta$ on $V$ such that $
\left|{{pL_{\lambda,v}}|_{M+\mathrm{i}\R}}^{} \right|\leq
\beta(v)$. By Corollary \ref{LineBound} we choose $N<0$ such for
any polynomial $p$ there exists a semi-norm $\beta'$ on $V$ such
that $ \left |{{pL_{\lambda,v}}|_{N+\mathrm{i}\R}}^{}
\right|\leq\beta'(v)$. Let $q$ be a polynomial such that the
multiset of poles of $1/q$ (with multiplicities) coincides with
the multiset of  poles of $L_{\pi}|_{[N,M]+\mathrm{i}\R}$. Here,
$[N,M]+\mathrm{i}\R$ denotes the strip $N\leq \re(s) \leq M$. It
is enough to show that the map $L'(\lambda)$ defined by
$L'(\lambda)(v):= L_{\lambda,v} q(\frac{1}{2})$ is continuous. Now
there exists a semi-norm $\alpha$ on $V$ such that, for any $v\in
V$, $$\left |qL_{\lambda,v}|_{M+\mathrm{i}\R} \right|\leq
\alpha(v) \text{ and } \left |qL_{\lambda,v}|_{N+\mathrm{i}\R}
\right|\leq \alpha(v).$$ By Theorem \ref{Bound}, for any $v \in
V$, there exists $\Delta$ such that $|qL_{\lambda,v}(s)| \leq
\alpha(v)$ if $s \in [N,M]+\mathrm{i}\R$ and $|\Im s|\geq \Delta$.
Now by maximal modulus principle $L'(\lambda)(v) \leq \alpha(v)$
for any $v \in V$.
\end{proof}

\begin{definition}
Let $(\pi,V)$ an irreducible admissible representation of $G$. We
define a map $$L: (V^*)^{HU,\Psi} \to (V^*)^{HG_1}$$ by
$$ L(\lambda)(v)=\frac{L_{\lambda,v}}{ L_\pi}(\frac{1}{2})$$
\end{definition}

Theorem \ref{IntShal} implies
\begin{proposition}
$L$ is a monomorphism.
\end{proposition}

Now we use the following theorem from \cite{AG_RJR}.
\begin{theorem}[see \cite{AG_RJR}, Theorem I] \label{UnLin}
The pair $(GL_{2n},GL_{n} \times GL_n)$ is a Gelfand pair. Namely,
$$\dim (V^*)^{\GL_n(F) \times \GL_n(F)} \leq 1.$$
\end{theorem}

\begin{corollary}
Theorem \ref{UShal} holds. Namely,
$$ \dim (V^*)^{HU,\Psi} \leq 1.$$
\end{corollary}


\end{document}